\documentclass[10pt]{amsart}
\usepackage{amsfonts}
\usepackage{amsmath}

\newtheorem{theorem}{Theorem}
\newtheorem{definition}{Definition}
\newtheorem{proposition}{Proposition}

\newtheorem{example}{Example}
\begin{document}
\title{Geometrical structures on the cotangent bundle}
\author{Liviu Popescu}
\maketitle
\begin{abstract}
In this paper we study the geometrical structures on the cotangent
bundle using the notions of adapted tangent structure and regular
vector fields. We prove that the dynamical covariant derivative on
$T^{*}M$ fix a nonlinear connection for a given
$\mathcal{J}$-regular vector field. Using the Legendre
transformation induced by a regular Hamiltonian, we show that a
semi-Hamiltonian vector field on $T^{*}M$ corresponds to a
semispray on $TM$ if and only if the nonlinear connection on $TM$
is just the canonical nonlinear connection induced by the regular
Lagrangian.
\end{abstract}
MSC2010: 53C05, 53C60, 70H05\\
Keywords: adapted tangent structure, regular vector field,
dynamical covariant derivative, Jacobi endomorphism, Hamiltonian
vector field.
\section{Introduction}

It is well known that many geometric structures on the cotangent
bundle $ T^{*}M$ of a differentiable manifold $M$ can be studied
using the same methods as in the case of tangent bundle $TM$. The
tangent bundle has a naturally defined integrable tangent
structure and together with a semispray (second order differential
equation vector field) induce a nonlinear connection \cite{Cr1,
Gr1}. However, in the case of the cotangent bundle we do not have
a canonical tangent structure or something similar to a semispray,
but there are several dual objects, as the adapted almost tangent
structure and regular vector fields \cite{Op}. Also, we have some
canonical geometric objects: Liouville-Hamilton vector field,
Liouville 1-form and symplectic structure \cite{Ab, Mi2, Ya}. The
existence of a regular Hamiltonian on $T^{*}M$ or a regular
Lagrangian on $TM$ permit us to transfer many geometrical
structures between these spaces using the Legendre transformation
\cite{Hr, Mi2, Po1}. In the paper \cite{Po2} the dynamical
covariant derivative on the cotangent bundle is introduced as a
tensor derivation and the metric nonlinear connections are
studied.\\ In this paper we investigate more geometrical
structures on the cotangent bundle and show how the dynamical
covariant derivative induced by a regular vector field and an
arbitrary nonlinear connection fix the nonlinear connection. The
paper is organized as follows. In the first section we present the
preliminary results on the cotangent bundle (see for instance
\cite{Id, Mi2, Mit, Op, Ya}). Also, we study the tension and the
strong torsion of the nonlinear connection and investigate the
homogeneous case. Using the Fr\"olisher-Nijenhuis bracket
\cite{Fr, Sz2} we study the Jacobi endomorphism on the cotangent
bundle. In the second section, using the notions of
$\mathcal{J}$-regular vector field and an arbitrary nonlinear
connection we introduce the dynamical covariant derivative on the
cotangent bundle and prove that the condition of compatibility
with the adapted tangent structure, that is $\nabla \mathcal{J}$
$=0,$ fix the nonlinear connection (see \cite{Bu3} for the tangent
bundle). In the section three we investigate the canonical
nonlinear connection induced by a $\mathcal{J}$-regular vector
field. We introduce the almost complex structure using the
Fr\"olisher-Nijenhuis bracket and prove that its dynamical
covariant derivative vanishes. These subjects can by found in the
case of tangent bundle in a lot of papers (see for instance
\cite{Bu2, Bu3, Cr2, Cr3, Gr2, Ku, Le, Ma, Sa1, Sa2, Sz1}).\\
In the last part, these geometrical structures are exemplified in
the case of Hamiltonian space. Finally, using the Legendre
transformation induced by a regular Hamiltonian, we prove that a
semispray on $TM$ corresponds to a semi-Hamiltonian vector field
on $T^{*}M$ if and only if the nonlinear connection on $TM$
determined by semispray is just the canonical nonlinear connection
induced by the regular Lagrangian.

\section{Preliminary structures on the cotangent bundle}

If $M$ be a differentiable, $n$-dimensional manifold and
$(T^{*}M,\tau ,M)$ the cotangent bundle then the local coordinates
on $\tau ^{-1}(U)$ are denoted $(x^i,p_i),$ $(i=\overline{1,n})$.
The natural basis on $T^{*}M$ is denoted $\left( \frac \partial
{\partial x^i},\frac \partial {\partial p_i}\right) $ and we
consider the following geometric objects
\begin{equation*}
C^{*}=p_i\frac \partial {\partial p_i},\quad \theta =p_idx^i,\quad
\omega =d\theta =dp_i\wedge dx^i,
\end{equation*}
where $(dx^i,dp_i)$ is the dual natural basis. The following
properties hold:\\
1$^{\circ }$ $C^{*}$ is a vertical vector field, globally defined
on $T^{*}M$, called the Liouville-Hamilton vector field.\\
2$^{\circ }$ The 1-form $\theta $ is globally defined on $T^{*}M$
and is called the Liouville 1-form.\\
3$^{\circ }$ The 2-form $\omega $ is a symplectic structure,
called canonical.\\
The Poisson bracket $\{\cdot ,\cdot \}$ on $T^{*}M$, is defined by
\begin{equation*}
\{f,g\}=\frac{\partial f}{\partial p_i}\frac{\partial g}{\partial
x^i}-\frac{
\partial g}{\partial p_i}\frac{\partial f}{\partial x^i},\quad \forall
f,g\in \mathcal{F}(T^{*}M).
\end{equation*}
In the following by a $d$-tensor field we mean a tensor field on
$T^{*}M$ whose components, under a change of coordinates on
$T^{*}M$, behave like the components of a tensor field on $M$. We
recall that if $L$ and $K$ are $(1,1)$-type tensor field,
Fr\"olisher-Nijenhuis bracket $[L,K]$ is the vector valued 2-form
\cite{Fr}
\begin{eqnarray*}
\lbrack L,K](X,Y) &=&[LX,KY]+[KX,LY]+(LK+KL)[X,Y] \\
&&\ \ \ -L[X,KY]-K[X,LY]-L[KX,Y]-K[LX,Y].
\end{eqnarray*}
and the Nijenhuis tensor of $L$ is given by
\begin{equation*}
N_L(X,Y)=\frac 12[L,L]=[LX,LY]+L^2[X,Y]-L[X,LY]-L[LX,Y].
\end{equation*}
For a vector field in $\mathcal{X}(M)$ and a $(1,1)$-type tensor
field $L$ on $M$ the Fr\"olisher-Nijenhuis bracket
$[X,L]=\mathcal{L}_XL$ is the $ (1,1) $-type tensor field on $M$
given by
\begin{equation*}
\mathcal{L}_XL=\mathcal{L}_X\circ L-L\circ \mathcal{L}_X,
\end{equation*}
where $\mathcal{L}_X$ is the usual Lie derivative.\\
On the cotangent bundle $T^{*}M$ there exists the integrable
vertical distribution $V_uT^{*}M$, $u\in T^{*}M$ generated locally
by the basis $ \left\{ \frac \partial {\partial p_i}\right\}
_{i=\overline{1,n}}$. A nonlinear connection $\mathcal{N}$ is a
horizontal distribution $H_uT^{*}M$
which is supplementary to the vertical distribution, that is $%
T_uT^{*}M=V_uT^{*}M\oplus H_uT^{*}M.$ If $\mathcal{N}$ is a
nonlinear connection then on the every domain of the local chart
$\tau ^{-1}(U)$, the adapted basis of the horizontal distribution
$HT^{*}M$ is \cite{Mi2}
\begin{equation*}
\frac \delta {\delta x^i}=\frac \partial {\partial
x^i}+\mathcal{N} _{ij}\frac \partial {\partial p_j},
\end{equation*}
where $\mathcal{N}_{ij}(x,p)$ are the coefficients of the
nonlinear connection. The dual adapted basis is
\begin{equation*}
\delta p_i=dp_i-\mathcal{N}_{ij}dx^j.
\end{equation*}
The system of vector fields $\left( \frac \delta {\delta
x^i},\frac \partial {\partial p_i}\right) $ defines the local
Berwald basis on $T^{*}M$. A nonlinear connection induces the
horizontal and vertical projectors given by
\begin{equation*}
h=\frac 12(Id+\mathcal{N}),\quad v=\frac 12(Id-\mathcal{N}),
\end{equation*}
and locally
\begin{equation*}
h=\frac \delta {\delta x^i}\otimes dx^i,\quad v=\frac \partial
{\partial p_i}\otimes \delta p_i.
\end{equation*}
The following properties hold
\begin{equation*}
h^2=h,\quad v^2=v,\quad hv=vh=0,\quad h+v=Id,\quad
h-v=\mathcal{N}.
\end{equation*}
We consider the nonlinear connection $\mathcal{N}$ and denote
$\tau _{ij}=\frac 12(\mathcal{N}_{ij}-\mathcal{N}_{ji}).$ The
nonlinear connection $\mathcal{N}$ on $T^{*}M$ is called symmetric
if
\begin{equation*}
\omega (hX,hY)=0,\quad X,Y\in \mathcal{X}(T^{*}M),
\end{equation*}
where $h$ is the horizontal projector induced by the nonlinear
connection. Locally, we obtain that the nonlinear connection is
symmetric if and only if $\tau _{ij}=0$, that is
$\mathcal{N}_{ij}=\mathcal{N}_{ji}.$ The following equations hold
\cite{Mi2}
\begin{equation}
\left[ \frac \delta {\delta x^i},\frac \delta {\delta x^j}\right]
=R_{kij}\frac \partial {\partial p_k},\quad \ \left[ \frac \delta
{\delta x^i},\frac \partial {\partial p_j}\right] =-\frac{\partial
\mathcal{N}_{ir}}{
\partial p_j}\frac \partial {\partial p_r},\quad \left[ \frac \partial
{\partial p_i},\frac \partial {\partial p_j}\right] =0,
\end{equation}
\begin{equation}
R_{ijk}=\frac{\delta \mathcal{N}_{jk}}{\delta x^i}-\frac{\delta
\mathcal{N} _{ik}}{\delta x^j}.
\end{equation}
The curvature of the nonlinear connection $\mathcal{N}$ on
$T^{*}M$ is given by $\Omega =-\frac 12[h,h],$ where $h$ is the
horizontal projector and $ \frac 12[h,h]$ is the Nijenhuis tensor
of $h$. In local coordinates we obtain
\begin{equation*}
\Omega =-\frac 12R_{ijk}dx^i\wedge dx^j\otimes \frac \partial
{\partial p_k},
\end{equation*}
where $R_{ijk}$ is given by (2) and is called the curvature
$d$-tensor of the nonlinear connection $\mathcal{N}$. The
curvature of a nonlinear connection is an obstruction to the
integrability of the horizontal distribution. Using (1), it
results that the horizontal distribution is integrable if and only
if the curvature vanishes. Also, from the Jacobi identity we
obtain
\begin{equation*}
\lbrack h,\Omega ]=0.
\end{equation*}

\begin{definition}
The tension of the nonlinear connection $\mathcal{N}$ is given by
\begin{equation*}
\mathbf{t}=\frac 12\mathcal{L}_{C^{*}}\mathcal{N}
\end{equation*}
\end{definition}

In local coordinates we obtain
\begin{equation*}
\mathbf{t}=\left( p_k\frac{\partial \mathcal{N}_{ij}}{\partial
p_k}-\mathcal{ N}_{ij}\right) dx^i\otimes \frac \partial {\partial
p_j}.
\end{equation*}
It results that the tension of the nonlinear connection vanishes
if and only if $\mathcal{N}$ is homogeneous of degree one with
respect to $p$.

\begin{definition}
An almost tangent structure on $T^{*}M$ is a morphism
$\mathcal{J}:\mathcal{X }(T^{*}M)$ $\rightarrow
\mathcal{X}(T^{*}M)$ of rank $n$ such that $\mathcal{ J}^2=0$. The
almost tangent structure is called adapted if \cite{Op}
\begin{equation*}
Im\mathcal{J}=Ker\mathcal{J}=VT^{*}M.
\end{equation*}
\end{definition}

The following properties hold
\begin{equation}
\mathcal{J}h=\mathcal{J},\quad h\mathcal{J}=0,\quad
\mathcal{J}v=0,\quad v \mathcal{J}=\mathcal{J}.
\end{equation}
Locally, an adapted almost tangent structure has the form
\begin{equation}
\mathcal{J}=t_{ij}dx^i\otimes \frac \partial {\partial p_j},
\end{equation}
where $t_{ij}(x,p)$ is a $d$-tensor field of rank $n.$ The adapted
almost tangent structure $\mathcal{J}$ is integrable if and only
if
\begin{equation}
\frac{\partial t^{ij}}{\partial p_k}=\frac{\partial
t^{kj}}{\partial p_i},
\end{equation}
where $t_{ij}t^{jk}=\delta _i^k.$ The adapted almost tangent
structure $ \mathcal{J}$ is called symmetric if
\begin{equation*}
\omega (\mathcal{J}X,Y)=\omega (\mathcal{J}Y,X).
\end{equation*}
Locally, this relation is equivalent with the symmetry of the
$d$-tensor $ t_{ij}(x,p)$.

\begin{example}
If $g$ is a pseudo-Riemannian metric on vertical subbundle
$VT^{*}M$, then there exists a unique adapted almost tangent
structure $\mathcal{J}$ on $ T^{*}M$ such that
\begin{equation}
g(\mathcal{J}X,\mathcal{J}Y)=-\omega (\mathcal{J}X,Y),\quad X,Y\in
\mathcal{X }(T^{*}M),
\end{equation}
and we say that $\mathcal{J}$ is induced by the metric $g$.
Locally, if we consider
\begin{equation*}
g(x,p)=g^{ij}dp_i\otimes dp_j,
\end{equation*}
then (6) implies that $t^{ij}=g^{ij}.$
\end{example}

\begin{definition}
An adapted almost tangent structure $\mathcal{J}$ is called
homogeneous if
\begin{equation}
\mathcal{L}_{C^{*}}\mathcal{J}=-\mathcal{J}.
\end{equation}
\end{definition}

\begin{proposition}
The adapted almost tangent structure is homogeneous if and only if
the local components $t_{ij}(x,p)$ are $0$-homogeneous with
respect to $p$.
\end{proposition}

\textbf{Proof}. By direct computation in local coordinates we
obtain
\begin{equation*}
\mathcal{L}_{C^{*}}\mathcal{J}=\left( p_k\frac{\partial
t_{ij}}{\partial p_k} -t_{ij}\right) dx^i\otimes \frac \partial
{\partial p_j},
\end{equation*}
and the equation (7) leads to $p_k\frac{\partial t_{ij}}{\partial
p_k}=0$, which ends the proof.\hfill\hbox{\rlap{$\sqcap$}$\sqcup$}

\begin{definition}
The torsion of a nonlinear connection $\mathcal{N}$ on $T^{*}M$ is
defined by $\mathcal{T}=[\mathcal{J},h]$, where $h$ is the
horizontal projector and $ [\mathcal{J},h]$ is the
Fr\"olicher-Nijenhuis bracket
\begin{equation*}
\begin{array}{c}
[\mathcal{J},h](X,Y)=[\mathcal{J}X,hY]+[hX,\mathcal{J}Y]+\mathcal{J}[X,Y]-
\mathcal{J}[X,hY]- \\
-\mathcal{J}[hX,Y]-h[X,\mathcal{J}Y]-h[\mathcal{J}X,Y].
\end{array}
\end{equation*}
\end{definition}

Locally, we consider
\begin{equation*}
\mathcal{T}=\frac 12\mathcal{T}_{ijk}(dx^i\wedge dx^j)\otimes
\frac \partial {\partial p_k},
\end{equation*}
and by straightforward computation, it results
\begin{equation}
\mathcal{T}_{ijk}=t_{is}\frac{\partial \mathcal{N}_{jk}}{\partial
p_s}-t_{js} \frac{\partial \mathcal{N}_{ik}}{\partial
p_s}+\frac{\delta t_{jk}}{\delta x^i}-\frac{\delta t_{ik}}{\delta
x^j}.
\end{equation}

\begin{proposition}
The following equations hold
\begin{equation*}
[\mathcal{J},\Omega ]=[h,\mathcal{T}],\quad [C^{*},\Omega
]=[h,\mathbf{t}].
\end{equation*}
\end{proposition}

\textbf{Proof}. We apply the Jacobi identity for the
Frolicher-Nijenhuis bracket
\begin{equation*}
\lbrack
\mathcal{J},[h,h]]+[h,[h,\mathcal{J}]]+[h,[\mathcal{J},h]]=0,
\end{equation*}
\begin{equation*}
-2[\mathcal{J},\Omega ]=-2[h,[h,\mathcal{J}]],
\end{equation*}
which yields $[\mathcal{J},\Omega ]=[h,\mathcal{T}].$ Also
\begin{equation*}
\lbrack C^{*},[h,h]]+[h,[h,C^{*}]]+[h,[C^{*},h]]=0,
\end{equation*}
\begin{equation*}
\lbrack C^{*},[h,h]]=-2[h,[C^{*},h]],
\end{equation*}
and it results $[C^{*},\Omega ]=[h,\mathbf{t}]$, because
$\mathbf{t}=\frac 12
\mathcal{L}_{C^{*}}\mathcal{N}=[C^{*},h].$\hfill\hbox{\rlap{$\sqcap$}$%
\sqcup$}\\
Let $\mathcal{J}$ be the adapted almost tangent structure on
$T^{*}M$. From \cite{Op} we set:

\begin{definition}
A vector field $\rho $ $\in \mathcal{X}(T^{*}M)$ is called
$\mathcal{J}$-regular if it satisfies the equation
\begin{equation}
\mathcal{J}[\rho ,\mathcal{J}X]=-\mathcal{J}X,\quad \forall X\in
\mathcal{X} (T^{*}M).
\end{equation}
\end{definition}

Locally, a vector field on $T^{*}M$ given in local coordinates by
\begin{equation*}
\rho =\xi ^i(x,p)\frac \partial {\partial x^i}+\chi _i(x,p)\frac
\partial {\partial p_i},
\end{equation*}
is $\mathcal{J}$-regular if and only if
\begin{equation}
t^{ij}=\frac{\partial \xi ^j}{\partial p_i},
\end{equation}
where $t_{ij}t^{jk}=\delta _i^k.$ We have to remark that, if the
equation $ \mathcal{J}[\rho ,\mathcal{J}X]=-\mathcal{J}X$ is
satisfied for any $X\in \mathcal{X}(T^{*}M),$ with the condition
$rank\left[ \frac{\partial \xi ^j}{
\partial p_i}\right] =n,$ then $\mathcal{J}$ is an integrable structure.

\begin{proposition}
If $\mathcal{J}$ is a homogeneous adapted tangent structure then a
vector field $\rho $ is $\mathcal{J}$-regular if and only if
\begin{equation*}
\mathcal{J}\rho =C^{*}.
\end{equation*}
\end{proposition}

\textbf{Proof}. Suppose that $\rho $ is $\mathcal{J}$-regular,
then $t^{ij}= \frac{\partial \xi ^j}{\partial p_i}$ and $t_{ij}$
is $0$-homogeneous, hence $\xi ^{j\text{ }}$ is $1$-homogeneous
with respect to $p$, therefore $p_i \frac{\partial \xi
^j}{\partial p_i}=\xi ^j$ that is $p_it^{ij}=\xi ^j$ which yields
$\mathcal{J}\rho =C^{*}$. Conversely, if $\mathcal{J}\rho =C^{*} $
then $\xi ^j=p_it^{ij}$ implies
\begin{equation*}
\frac{\partial \xi ^j}{\partial p_k}=\frac{\partial
t^{ij}}{\partial p_k} p_i+t^{kj}=\frac{\partial t^{kj}}{\partial
p_i}p_i+t^{kj}=t^{kj}
\end{equation*}
so $\rho $ is
$\mathcal{J}$-regular.\hfill\hbox{\rlap{$\sqcap$}$\sqcup$}

\begin{definition}
The strong torsion of the nonlinear connection has the form
\begin{equation}
\Bbb{T}=i_\rho \mathcal{T}-\mathbf{t}.
\end{equation}
where $i_\rho $ is the interior product with respect to the
$\mathcal{J}$-regular vector field.
\end{definition}

In local coordinates we obtain that
\begin{equation}
\Bbb{T}=\left( \xi
^i\mathcal{T}_{ijk}+\mathcal{N}_{jk}-p_s\frac{\partial
\mathcal{N}_{jk}}{\partial p_s}\right) dx^j\otimes \frac \partial
{\partial p_k},
\end{equation}
with $\mathcal{T}_{ijk}$ from (8). Using the relations (11) and
(12) we obtain:

\begin{proposition}
The strong torsion of a nonlinear connection vanishes if the weak
torsion and the tension vanish. Conversely is true if the
nonlinear connection is homogeneous of degree one.
\end{proposition}

\begin{definition}
For a $\mathcal{J}$-regular vector field $\rho $ and an arbitrary
nonlinear connection $\mathcal{N}$ with induced $(h,v)\,$
projectors, we consider the vertically valued $(1,1)$-type tensor
field on $T^{*}M\backslash \{0\}$ given by
\begin{equation*}
\Phi =v\circ \mathcal{L}_\rho h,
\end{equation*}
which will be called the Jacobi endomorphism.
\end{definition}

It results that
\begin{equation*}
\Phi =v\circ \mathcal{L}_\rho h=-v\circ \mathcal{L}_\rho v=v\circ
(\mathcal{L }_\rho \circ h-h\circ \mathcal{L}_\rho )=v\circ
\mathcal{L}_\rho \circ h.
\end{equation*}
In local coordinates we obtain
\begin{equation*}
\mathcal{L}_\rho \frac \delta {\delta x^j}=-\frac{\delta \xi
^i}{\delta x^j} \frac \delta {\delta x^i}+R_{jk}\frac \partial
{\partial p_k},\quad \mathcal{ L}_\rho \frac \partial {\partial
p_j}=-t^{ji}\frac \delta {\delta x^i}+\left(
t^{ji}\mathcal{N}_{ik}-\frac{\partial \chi _k}{\partial p_j}
\right) \frac \partial {\partial p_k}.
\end{equation*}
Locally, the Jacobi endomorphism has the form
\begin{equation}
\Phi =\mathcal{R}_{ij}dx^i\otimes \frac \partial {\partial
p_j},\quad \mathcal{R}_{jk}=\frac{\delta \xi ^i}{\delta
x^j}\mathcal{N}_{ik}-\frac{ \delta \chi _k}{\delta x^j}+\rho
(\mathcal{N}_{jk}).
\end{equation}

\begin{proposition}
The following result holds
\begin{equation*}
\Phi =i_\rho \Omega +v\circ \mathcal{L}_{v\rho }h.
\end{equation*}
\end{proposition}

\textbf{Proof}. Indeed, $\Phi (X)=v[\rho ,hX]=v[h\rho ,hX]+v[v\rho
,hX]$ and $\Omega (\rho ,X)=v[h\rho ,hX],$ which yields $\Phi
(X)=\Omega (\rho ,X)+v[v\rho ,hX].$
\hfill\hbox{\rlap{$\sqcap$}$\sqcup$}

\begin{proposition}
If $\rho $ is a horizontal $\mathcal{J}$-regular vector field then
\begin{equation*}
\Phi =i_\rho \Omega .
\end{equation*}
\end{proposition}

\textbf{Proof}. We have $\rho =h\rho $ and $v\rho =0$ which yields
$\Phi =i_\rho \Omega $ and locally it results
\begin{equation*}
\rho =\xi ^i\frac \delta {\delta x^i},\quad \chi _i=\xi
^k\mathcal{N} _{ki},\quad \mathcal{R}_{ij}=R_{kij}\xi ^k,
\end{equation*}
which show us the relation between the Jacobi endomorphism given
by (13) and curvature tensor from (2).
\hfill\hbox{\rlap{$\sqcap$}$\sqcup$}\\
We recall that the pair $(T^{*}M,\omega )$ is a symplectic
manifold. Every differentiable function $f:T^{*}M\rightarrow R$
determines a unique vector field $X_f$ on $T^{*}M$ called the
Hamiltonian vector field such that $ i_{X_f}\omega =-df$. It
results that $\mathcal{L}_{X_f}\omega =0$ and considering a
$\mathcal{J}$-regular Hamiltonian vector field $\rho =\xi ^i\frac
\partial {\partial x^i}+\chi _i\frac \partial {\partial p_i}$ this
condition is equivalent with \cite{Op}
\begin{equation}
a)\ \frac{\partial \xi ^j}{\partial p_i}=\frac{\partial \xi
^i}{\partial p_j} ,\quad b)\ \frac{\partial \chi _i}{\partial
p_j}=-\frac{\partial \xi ^j}{
\partial x^i},\quad c)\ \frac{\partial \chi _i}{\partial x^j}=\frac{\partial
\chi _j}{\partial x^i}.
\end{equation}
A vector field $X\in \mathcal{X}($ $T^{*}M)$ is a semi-Hamiltonian
vector field if it is $\mathcal{J}$-regular and
$\mathcal{L}_X\omega $ is a semibasic 2-form on $T^{*}M$, i.e.
satisfies the equation
\begin{equation*}
i_\nu (\mathcal{L}_X\omega )=0,\quad \forall \nu \in \Gamma
(VT^{*}M).
\end{equation*}
By direct computation, it results that in the case of
semi-Hamiltonian vector field, only the conditions (14) $a)$ and
$b)$ are satisfied.

\section{Dynamical covariant derivative on the cotangent bundle}

The notion of dynamical covariant derivative was introduced for
the first time in the case of tangent bundle by Cari\~nena and
Martinez \cite{Ca} as a derivation of degree 0 along the tangent
bundle projection (see also \cite {Bu2, Bu3, Cr2, Ma, Ma1, Sz2}).
An extensive study and discussions about the dynamical covariant
derivative which is associated to a second order vector field
(semispray) on $TM$ can be found in \cite{Sa2}. We introduce the
dynamical covariant derivative on the cotangent bundle which is
associated to a $\mathcal{J}$-regular vector field $\rho $ and an
arbitrary nonlinear connection $\mathcal{N}$, that is not fixed
yet. We determine the nonlinear connection by requiring the
compatibility of the dynamical covariant derivative with the
adapted tangent structure. Using \cite {Bu3, Po2} we set:

\begin{definition}
A map $\nabla :\mathcal{T}(T^{*}M\backslash \{0\})\rightarrow
\mathcal{T}
(T^{*}M\backslash \{0\})$ is said to be a tensor derivation on $%
T^{*}M\backslash \{0\}$ if the following conditions are
satisfied:\\
i) $\nabla $ is $\Bbb{R}$-linear,\\
ii) $\nabla $ is type preserving, i.e. $\nabla (\mathcal{T}%
_s^r(T^{*}M\backslash \{0\})\subset
\mathcal{T}_s^r(T^{*}M\backslash \{0\})$ , for each $(r,s)\in
\Bbb{N}\times \Bbb{N.},$\\
iii) $\nabla $ obeys the Leibnitz rule $\nabla (T\otimes S)=\nabla
T\otimes S+T\otimes \nabla S,$\\
iv) $\nabla \,$commutes with any contractions.
\end{definition}

For a $\mathcal{J}$-regular vector field $\rho $ and an arbitrary
nonlinear connection $\mathcal{N}$ with induced $(h,v)\,$
projectors, we consider the map $\nabla
:\mathcal{X}(T^{*}M\backslash \{0\})\rightarrow \mathcal{X}
(T^{*}M\backslash \{0\})$ given by
\begin{equation}
\nabla =h\circ \mathcal{L}_\rho \circ h+v\circ \mathcal{L}_\rho
\circ v,
\end{equation}
which is called the dynamical covariant derivative with respect to
$\rho $ and the nonlinear connection $\mathcal{N}$. We have to
remark that in this moment we will not yet assume any dependence
between $\rho $ and the nonlinear connection $\mathcal{N}$ with
induced $(h,v)$ projectors. By setting $\nabla f=\rho (f),$ for
$f\in C^\infty (T^{*}M\backslash \{0\})$ using the Leibnitz rule
and the requirement that $\bigtriangledown $ commutes with any
contraction, we can extend the action of $\bigtriangledown $ to
arbitrary tensor fields on $T^{*}M\backslash \{0\}$. For a 1-form
on $ T^{*}M\backslash \{0\}$ the dynamical covariant derivative is
given by
\begin{equation*}
(\nabla \varphi )(X)=\rho (\varphi )(X)-\varphi (\nabla X).
\end{equation*}
For a $(1,1)$-type tensor field $T$ on $T^{*}M\backslash \{0\}$
the dynamical covariant derivative has the form
\begin{equation*}
\nabla T=\nabla \circ T-T\circ \nabla .
\end{equation*}
By direct computation we obtain $\nabla h=\nabla v=0$ and
\begin{equation}
\nabla \frac \delta {\delta x^j}=-\frac{\delta \xi ^i}{\delta
x^j}\frac \delta {\delta x^i},\quad \nabla dx^j=\frac{\delta \xi
^j}{\delta x^i}dx^i,
\end{equation}
\begin{equation}
\nabla \frac \partial {\partial p_j}=\left(
t^{ji}\mathcal{N}_{ik}-\frac{
\partial \chi _k}{\partial p_j}\right) \frac \partial {\partial p_k},\quad
\nabla \delta p_j=-\left( t^{ki}\mathcal{N}_{ij}-\frac{\partial
\chi _j}{
\partial p_k}\right) \delta p_k.
\end{equation}
Next proposition proves the compatibility conditions between some
of the geometric structures on the cotangent bundle.

\begin{proposition}
The following results hold:
\begin{equation}
h\circ \mathcal{L}_\rho \circ \mathcal{J}=-h,\quad
\mathcal{J}\circ \mathcal{ L}_\rho \circ v=-v,
\end{equation}
\begin{equation}
\nabla \mathcal{J}=\mathcal{L}_\rho \mathcal{J}+h-v,\quad \nabla
\mathcal{J} =\left( \rho (t_{ij})+t_{kj}\frac{\partial \xi
^k}{\partial x^i}-t_{ik}\frac{
\partial \chi _j}{\partial p_k}+2\mathcal{N}_{ij}\right) dx^i\otimes \frac
\partial {\partial p_j}.
\end{equation}
\end{proposition}

\textbf{Proof}. From (9) we get
\begin{equation*}
\mathcal{J}[\rho ,\mathcal{J}X]=-\mathcal{J}X\Rightarrow
\mathcal{J}\left( [\rho ,\mathcal{J}X]+X\right) =0\Rightarrow
[\rho ,\mathcal{J}X]+X\in VT^{*}M
\end{equation*}
\begin{equation*}
h\left( [\rho ,\mathcal{J}X]+X\right) =0\Rightarrow h[\rho
,\mathcal{J} X]=-hX\Leftrightarrow h\circ \mathcal{L}_\rho \circ
\mathcal{J}=-h.
\end{equation*}
Also, in $\mathcal{J}[\rho ,\mathcal{J}X]+\mathcal{J}X=0$
considering $ \mathcal{J}X=vZ$ it results $\mathcal{J}[\rho
,vZ]=-vZ\Leftrightarrow \mathcal{J}\circ \mathcal{L}_\rho \circ
v=-v.$ In the following we have
\begin{eqnarray*}
\nabla \circ \mathcal{J} &=&h\circ \mathcal{L}_\rho \circ h\circ
\mathcal{J} +v\circ \mathcal{L}_\rho \circ v\circ
\mathcal{J}=v\circ \mathcal{L}_\rho
\circ \mathcal{J}= \\
&=&(Id-h)\circ \mathcal{L}_\rho \circ \mathcal{J}=\mathcal{L}_\rho
\circ \mathcal{J}-h\circ \mathcal{L}_\rho \circ
\mathcal{J}=\mathcal{L}_\rho \circ \mathcal{J}+h.
\end{eqnarray*}
But, on the other hand
\begin{equation*}
\mathcal{J}\circ \nabla =\mathcal{J}\circ \mathcal{L}_\rho \circ
h=\mathcal{J }\circ \mathcal{L}_\rho \circ (Id-v)=\mathcal{J}\circ
\mathcal{L}_\rho - \mathcal{J}\circ \mathcal{L}_\rho \circ
v=\mathcal{J}\circ \mathcal{L}_\rho +v.
\end{equation*}
and it results
\begin{equation*}
\nabla \circ \mathcal{J}-\mathcal{J}\circ \nabla =\mathcal{L}_\rho
\circ \mathcal{J}+h-\mathcal{J}\circ \mathcal{L}_\rho
-v\Rightarrow \nabla \mathcal{J}=\mathcal{L}_\rho \mathcal{J}+h-v.
\end{equation*}
For the last relation, we have

\begin{eqnarray*}
\nabla \mathcal{J} &=&\nabla \left( t_{ij}dx^i\otimes \frac
\partial {\partial p_j}\right) =\nabla t_{ij}dx^i\otimes \frac
\partial {\partial
p_j}+t_{ij}dx^i\otimes \nabla \frac \partial {\partial p_j} \\
\  &=&(\rho (t_{ij})dx^i+t_{ij}\nabla dx^i)\otimes \frac \partial
{\partial
p_j}+t_{ij}dx^i\otimes \nabla \frac \partial {\partial p_j} \\
\  &=&\rho (t_{ij})dx^i\otimes \frac \partial {\partial
p_j}+t_{ij}\frac{\delta \xi ^i}{\delta x^k}dx^k\otimes \frac
\partial {\partial p_j}+t_{ij}dx^i\otimes \left(
t^{sj}\mathcal{N}_{sk}-\frac{\partial \chi _k}{
\partial p_j}\right) \frac \partial {\partial p_k} \\
\  &=&\left( \rho (t_{ij})+t_{kj}\frac{\delta \xi ^k}{\delta
x^i}+\mathcal{N} _{ij}-t_{ik}\frac{\partial \chi _j}{\partial
p_k}\right) dx^i\otimes \frac
\partial {\partial p_j} \\
\  &=&\left( \rho (t_{ij})+t_{kj}\left( \frac{\partial \xi
^k}{\partial x^i}+
\mathcal{N}_{is}\frac{\partial \xi ^k}{\partial p_s}\right) +\mathcal{N}%
_{ij}-t_{ik}\frac{\partial \chi _j}{\partial p_k}\right)
dx^i\otimes \frac
\partial {\partial p_j} \\
\  &=&\left( \rho (t_{ij})+t_{kj}\frac{\partial \xi ^k}{\partial x^i}%
+t_{kj}t^{sk}\mathcal{N}_{is}+\mathcal{N}_{ij}-t_{ik}\frac{\partial
\chi _j}{
\partial p_k}\right) dx^i\otimes \frac \partial {\partial p_j} \\
\  &=&\left( \rho (t_{ij})+t_{kj}\frac{\partial \xi ^k}{\partial
x^i}-t_{ik} \frac{\partial \chi _j}{\partial
p_k}+2\mathcal{N}_{ij}\right) dx^i\otimes \frac \partial {\partial
p_j}.
\end{eqnarray*}
\hfill\hbox{\rlap{$\sqcap$}$\sqcup$}.\\
We will use the previous results in order to fix a nonlinear
connection. We note that the formula $\nabla
\mathcal{J}=\mathcal{L}_\rho \mathcal{J}+ \mathcal{N}$ represents
the difference between an arbitrary nonlinear connection
$\mathcal{N}$, which is not fixed yet, and the canonical nonlinear
connection associated to a $\mathcal{J}$-regular vector field $
\rho,$ that will be fixed in the next theorem.

\begin{theorem}
For a $\mathcal{J}$-regular vector field $\rho $, an arbitrary
nonlinear connection $\mathcal{N}$ and the dynamical covariant
derivative $\nabla $ with respect to $\rho $ and $\mathcal{N}$
the following conditions are equivalent:\\
$i)$ $\nabla \mathcal{J}=0,$\\
$ii)$ $\mathcal{L}_\rho \mathcal{J}+h-v=0,$\\
$iii)$ $\mathcal{N}_{ij}=\frac 12\left( t_{ik}\frac{\partial \chi
_j}{
\partial p_k}-t_{kj}\frac{\partial \xi ^k}{\partial x^i}-\rho
(t_{ij})\right) .$
\end{theorem}

\textbf{Proof}. The proof follows from the relation (19).\hfill
\hbox{\rlap{$\sqcap$}$\sqcup$}\\
This theorem shows that the compatibility condition $\nabla
\mathcal{J}=0$ of the dynamical covariant derivative with the
adapted tangent structure $ \mathcal{J}$, fix the nonlinear
connection $\mathcal{N}=-\mathcal{L}_\rho \mathcal{J}$. In the
following we will deal with this nonlinear connection induced by
the $\mathcal{J}$-regular vector field $\rho $, called the
canonical nonlinear connection.

\subsection{The canonical nonlinear connection induced by a $
\mathcal{J}$-regular vector field}

Given an adapted tangent structure $\mathcal{J}$ and a
$\mathcal{J}$-regular vector field $\rho $, then the compatibility
condition $\nabla \mathcal{J}=0$ fix the canonical nonlinear
connection
\begin{equation}
\mathcal{N}=-\mathcal{L}_\rho \mathcal{J},
\end{equation}
with local coefficients given by \cite{Op}
\begin{equation}
\mathcal{N}_{ij}=\frac 12\left( t_{ik}\frac{\partial \chi
_j}{\partial p_k} -t_{kj}\frac{\partial \xi ^k}{\partial x^i}-\rho
(t_{ij})\right) .
\end{equation}

\begin{definition}
The almost complex structure has the form
\begin{equation*}
\Bbb{F}=h\circ \mathcal{L}_\rho h-\mathcal{J}.
\end{equation*}
\end{definition}

\begin{proposition}
The following results hold:
\begin{eqnarray*}
\Bbb{F}^2=-Id,\quad \Bbb{F}\circ \mathcal{J}=h,\quad
\mathcal{J}\circ \Bbb{F} =v,\quad v\circ \Bbb{F}=\Bbb{F}\circ
h=-\mathcal{J}, \\ \quad h\circ \Bbb{F}=\Bbb{F}\circ
v=\Bbb{F}+\mathcal{J},\quad \mathcal{N} \circ
\Bbb{F}=\Bbb{F}+2\mathcal{J},\quad \Phi =\mathcal{L}_\rho
h-\Bbb{F}- \mathcal{J}.
\end{eqnarray*}
\end{proposition}

\textbf{Proof}. It results $\Bbb{F}=h\circ \mathcal{L}_\rho \circ
h-h\circ \mathcal{L}_\rho -\mathcal{J}$ $=h\circ \mathcal{L}_\rho
\circ (h-Id)-J=-h\circ \mathcal{L}_\rho \circ v-\mathcal{J}$.\\
$\Bbb{F}^2=\left( -h\circ \mathcal{L}_\rho \circ
v-\mathcal{J}\right) \circ \left( -h\circ \mathcal{L}_\rho \circ
v-\mathcal{J}\right) =h\circ \mathcal{L }_\rho \circ v\circ h\circ
\mathcal{L}_\rho \circ v+h\circ \mathcal{L}_\rho \circ v\circ
\mathcal{J}+$\\
$+\mathcal{J}\circ h\circ \mathcal{L}_\rho \circ
v+\mathcal{J}^2=h\circ \mathcal{L}_\rho \circ
\mathcal{J}+\mathcal{J}\circ \mathcal{L}_\rho \circ v$
$=-h-v=-Id.$\\
$\Bbb{F}\circ \mathcal{J}=\left( -h\circ \mathcal{L}_\rho \circ
v-\mathcal{J} \right) \circ \mathcal{J}=$ $-h\circ
\mathcal{L}_\rho \circ v\circ \mathcal{J }-\mathcal{J}^2$
$=-h\circ \mathcal{L}_\rho \circ \mathcal{J}=h.$\\
In the same way, using (18) the other relations can be
proved.\hfill \hbox{\rlap{$\sqcap$}$\sqcup$}\\
In local coordinates we have
\begin{equation*}
\Bbb{F}=t^{ij}\frac \delta {\delta x^i}\otimes \delta
p_j-t_{ij}\frac
\partial {\partial p_i}\otimes dx^j.
\end{equation*}
For a $\mathcal{J}$-regular vector field $\rho $, and induced
nonlinear connection $\mathcal{N}$ we consider the
$\Bbb{R}$-linear map $\nabla _0: \mathcal{X}(T^{*}M\backslash
\{0\})\rightarrow \mathcal{X}(T^{*}M\backslash \{0\})$ given by
\begin{equation*}
\nabla _0X=h[\rho ,hX]+v[\rho ,vX],\quad \forall X\in \mathcal{X}
(T^{*}M\backslash \{0\}).
\end{equation*}
It results that
\begin{equation*}
\nabla _0(fX)=\rho (f)X+f\nabla _0X,\quad \forall f\in C^\infty
(T^{*}M\backslash \{0\}),\ X\in \mathcal{X}(T^{*}M\backslash
\{0\}).
\end{equation*}
Any tensor derivation on $T^{*}M\backslash \{0\}$ is completely
determined by its actions on smooth functions and vector fields on
$T^{*}M\backslash \{0\}$ (see \cite{Sz2} generalized Willmore's
theorem). Therefore there exists a unique tensor derivation
$\nabla $ on $T^{*}M\backslash \{0\}$ such that
\begin{equation*}
\nabla \mid _{C^\infty (T^{*}M\backslash \{0\})}=\rho ,\quad
\nabla \mid _{ \mathcal{X}(T^{*}M\backslash \{0\})}=\nabla _0.
\end{equation*}
We will call the tensor derivation $\nabla $, the
\textit{dynamical covariant derivative} induced by the
$\mathcal{J}$-regular vector field $ \rho.$\\
To express the action of $\nabla $, we note that for $\rho =\xi
^i\frac
\partial {\partial x^i}+\chi _i\frac \partial {\partial p_i}$ it results
\begin{equation}
\left[ \rho ,\frac \partial {\partial p_j}\right] =-t^{ij}\frac
\delta {\delta x^i}+\left( t^{ij}\mathcal{N}_{ik}-\frac{\partial
\chi _k}{\partial p_j}\right) \frac \partial {\partial p_k},
\end{equation}
\begin{equation}
\left[ \rho ,\frac \delta {\delta x^j}\right] =-\frac{\delta \xi
^i}{\delta x^j}\frac \delta {\delta x^i}+\mathcal{R}_{jk}\frac
\partial {\partial p_k},
\end{equation}
where $\mathcal{R}_{jk}$ are the components of the Jacobi
endomorphism.\\
The action of the dynamical covariant derivative on the Berwald
basis has the form
\begin{equation*}
\nabla \frac \delta {\delta x^j}=h\left[ \rho ,\frac \delta
{\delta x^j}\right] =-\frac{\delta \xi ^i}{\delta x^j}\frac \delta
{\delta x^i},
\end{equation*}
\begin{equation*}
\nabla \frac \partial {\partial p_j}=v\left[ \rho ,\frac \partial
{\partial p_j}\right] =\left(
t^{ij}\mathcal{N}_{ik}-\frac{\partial \chi _k}{\partial
p_j}\right) \frac \partial {\partial p_k}.
\end{equation*}

\begin{proposition}
The dynamical covariant derivative is given by
\begin{equation*}
\nabla =\mathcal{L}_\rho +\Bbb{F}+\mathcal{J}-\Phi .
\end{equation*}
\end{proposition}

\textbf{Proof}. Using (15) we obtain
\begin{eqnarray*}
\nabla &=&h\circ \mathcal{L}_\rho \circ h+v\circ \mathcal{L}_\rho \circ v\\
&=&h\circ \left( \mathcal{L}_\rho h+h\circ \mathcal{L}_\rho
\right) +v\circ
\left( \mathcal{L}_\rho v+v\circ \mathcal{L}_\rho \right) \\
&=&h\circ \mathcal{L}_\rho h+v\circ \mathcal{L}_\rho v+(h+v)\circ
\mathcal{L}
_\rho \\
&=&\mathcal{L}_\rho +\Bbb{F}+\mathcal{J}-\Phi .
\end{eqnarray*}

\begin{proposition}
The following properties hold:
\begin{equation*}
\nabla \mathcal{J}=0,\ \nabla \Bbb{F}=0.
\end{equation*}
\end{proposition}

\textbf{Proof}. From (19) it results $\nabla \mathcal{J}=0$. From
$\Bbb{F} =-h\circ \mathcal{L}_\rho \circ v-\mathcal{J}$ and
$\nabla \Bbb{F}=\nabla \circ \Bbb{F}-\Bbb{F}\circ \nabla $ we
obtain

\begin{equation*}
\begin{array}{l}
\nabla \Bbb{F}=(h\circ \mathcal{L}_\rho \circ h+v\circ
\mathcal{L}_\rho
\circ v)\circ (-h\circ \mathcal{L}_\rho \circ v) \\
\qquad -(-h\circ \mathcal{L}_\rho \circ v)\circ (h\circ
\mathcal{L}_\rho
\circ h+v\circ \mathcal{L}_\rho \circ v) \\
\qquad =-h\circ \mathcal{L}_\rho \circ h\circ \mathcal{L}_\rho
\circ
v+h\circ \mathcal{L}_\rho \circ v\circ \mathcal{L}_\rho \circ v \\
\qquad =h\circ \mathcal{L}_\rho \circ (v-h)\circ \mathcal{L}_\rho
\circ v=h\circ \mathcal{L}_\rho \circ \mathcal{L}_\rho
\mathcal{J}\circ \mathcal{L}
_\rho \circ v \\
\qquad =h\circ \mathcal{L}_\rho \circ (\mathcal{L}_\rho \circ
\mathcal{J}-
\mathcal{J}\circ \mathcal{L}_\rho )\circ \mathcal{L}_\rho \circ v \\
\qquad =h\circ \mathcal{L}_\rho \circ \mathcal{L}_\rho \circ
(\mathcal{J} \circ \mathcal{L}_\rho \circ v)-(h\circ
\mathcal{L}_\rho \circ \mathcal{J}
)\circ \mathcal{L}_\rho \circ \mathcal{L}_\rho \circ v \\
\qquad =-h\circ \mathcal{L}_\rho \circ \mathcal{L}_\rho \circ
v+h\circ \mathcal{L}_\rho \circ \mathcal{L}_\rho \circ v=0.
\end{array}
\end{equation*}
using relations (18).\hfill \hbox{\rlap{$\sqcap$}$\sqcup$}\\
We have to remark that if $\rho $ is a horizontal
$\mathcal{J}$-regular vector field then $\nabla \rho =0$. Indeed
$\rho =h\rho $, $v\rho =0$ and it results $\nabla \rho =h\circ
\mathcal{L}_\rho \circ h\rho +v\circ \mathcal{L} _\rho \circ v\rho
=h\circ \mathcal{L}_\rho \circ \rho =0.$

\subsubsection{Application to Hamiltonian case}

A Hamilton space \cite{Mi2} is a pair $(M,H)$ where $M$ is a
differentiable, $n-$dimensional manifolds and $H$ is a function on
$T^{*}M$ with the properties:\\
1$^{\circ }$ $H:(x,p)\in T^{*}M\rightarrow H(x,p)\in \Bbb{R}$ is
differentiable on $T^{*}M$ and continue on the null section of the
projection $\tau :T^{*}M\rightarrow M$.\\
2$^{\circ }$ The Hessian of $H$ with respect to $p_i$ is
nondegenerate
\begin{equation}
g^{ij}=\frac{\partial ^2H}{\partial p_i\partial p_j},\ \
rank\left\| g^{ij}(x,p)\right\| =n\text{ on
}\widetilde{T^{*}M}=T^{*}M\backslash \{0\},
\end{equation}
3$^{\circ }$ $d$-tensor field $g^{ij}(x,p)$ has constant signature
on $ \widetilde{T^{*}M}.$\\
The Hamiltonian $H$ on $T^{*}M$ induces a pseudo-Riemannian metric
$g_{ij}$ with $g_{ij}g^{jk}=\delta _i^k$ and $g^{jk}$ given by
(24) on $VT^{*}M$. It induces a unique adapted tangent structure,
denoted
\begin{equation*}
\mathcal{J}_H=g_{ij}dx^i\otimes \frac \partial {\partial p_j},
\end{equation*}
such that (5) is satisfied. A $\mathcal{J}$-regular vector field
induced by the regular Hamiltonian $H$ is given by
\begin{equation*}
\rho _H=\frac{\partial H}{\partial p_i}\frac \partial {\partial
x^i}+\chi _i\frac \partial {\partial p_i}.
\end{equation*}
There exists a unique Hamiltonian vector field $\rho _H\in
\mathcal{X} (T^{*}M)$ such that
\begin{equation*}
i_{\rho _H}\omega =-dH,
\end{equation*}
given by
\begin{equation}
\rho _H=\frac{\partial H}{\partial p_i}\frac \partial {\partial
x^i}-\frac{
\partial H}{\partial x^i}\frac \partial {\partial p_i}.
\end{equation}
The nonlinear connection
\begin{equation*}
\mathcal{N}=-\mathcal{L}_{\rho _H}\mathcal{J}_H,
\end{equation*}
has the coefficients given by \cite{Mi2, Po1}
\begin{equation}
\mathcal{N}_{ij}=\frac 12\left( \{g_{ij},H\}-\left(
g_{ik}\frac{\partial ^2H }{\partial p_k\partial
x^j}+g_{jk}\frac{\partial ^2H}{\partial p_k\partial x^i}\right)
\right) ,
\end{equation}
and is called the \textit{canonical nonlinear connection} of the
Hamilton space $(M,H)$ which is a metric nonlinear connection
\cite{Po2}.\\
Let us consider the regular Hamiltonian $H(x,p)$ on $T^{*}M$ which
induces a local diffeomorphism $\Psi :T^{*}M\rightarrow TM$ given
by
\begin{equation}
x^i,\quad y^i=\xi ^i(x,p)=\frac{\partial H}{\partial p_i},
\end{equation}
and $\Psi ^{-1}$ has the following components
\begin{equation}
x^i,\quad p_i=\zeta _i(x,y)=\frac{\partial L}{\partial y^i},
\end{equation}
where
\begin{equation}
L(x,y)=\zeta _iy^i-H(x,p),
\end{equation}
is the Legendre transformation. From the condition for $\Psi
^{-1}$ to be the inverse of $\Psi $ we obtain \cite{Op}
\begin{equation*}
\frac{\partial \zeta _i}{\partial y^j}\circ \Psi =g_{ij},\quad
\frac{
\partial \zeta _i}{\partial x^j}\circ \Psi =-g_{ik}\frac{\partial \xi ^k}{
\partial x^j},\quad g^{ij}=\frac{\partial \xi ^j}{\partial p_i}=\frac{
\partial ^2H}{\partial p_i\partial p_j},\quad g_{ij}g^{jk}=\delta _i^k.
\end{equation*}
\begin{equation*}
\Psi _{*}\left( \frac \partial {\partial p_i}\right) =(g^{ik}\circ
\Psi ^{-1})\frac \partial {\partial y^k},\quad \Psi _{*}\left(
\frac \partial {\partial x^i}\right) =\frac \partial {\partial
x^i}+\left( \frac{\partial \xi ^k}{\partial x^i}\circ \Psi
^{-1}\right) \frac \partial {\partial y^k},
\end{equation*}
\begin{equation}
\Psi _{*}^{-1}\left( \frac \partial {\partial y^i}\right)
=g_{ki}\frac
\partial {\partial p_k},\quad \Psi _{*}^{-1}\left( \frac \partial {\partial
x^i}\right) =\frac \partial {\partial x^i}-g_{kh}\frac{\partial
\xi ^h}{
\partial x^i}\frac \partial {\partial p_k},
\end{equation}
where $\Psi _{*}$ is the tangent application of $\Psi $.\\
We know \cite{Cr1, Gr1} that if $S=y^i\frac \partial {\partial
x^i}+S^i\frac
\partial {\partial y^i}$ is a semispray and $J=dx^i\otimes \frac \partial
{\partial y^i}$ is the canonical tangent structure on $TM,$ then
the automorphism
\begin{equation*}
N=-\mathcal{L}_SJ,
\end{equation*}
is a nonlinear connection on $TM$ with the coefficients given by
\begin{equation}
N_j^i(x,y)=-\frac 12\frac{\partial S^i}{\partial y^j}.
\end{equation}
For a regular Lagrangian $L$ on $TM$ there exists the canonical
nonlinear connection with the coefficients (31) where
\begin{equation}
S^i=g^{ij}\left( \frac{\partial L}{\partial x^j}-\frac{\partial
^2L}{
\partial x^k\partial y^j}y^k\right) ,
\end{equation}
which is a metric nonlinear connection \cite{Bu1}. The next result
shows the relation between the canonical nonlinear connection
induces by the regular Lagrangian on $TM$ and the semi-Hamiltonian
vector field on \thinspace $ T^{*}M$, via Legendre transformation.
Let us consider a semispray $S$ on $TM$ and $\Psi ^{-1}$ the local
diffeomorphism given by (28).

\begin{theorem}
The vector field $\rho =\Psi _{*}^{-1}S$ is a semi-Hamiltonian
$\mathcal{J}$-regular vector field on $T^{*}M$ if and only if the
nonlinear connection on $TM$ induces by semispray is the canonical
nonlinear connection induced by the regular Lagrangian.
\end{theorem}

\textbf{Proof. }We consider a semispray $S=y^i\frac \partial
{\partial x^i}+S^i\frac \partial {\partial y^i}$ on $TM$ and from
(30) it results
\begin{equation*}
\Psi _{*}^{-1}S=\xi ^i\frac \partial {\partial q^i}+\left( -\xi
^ig_{kj} \frac{\partial \xi ^j}{\partial x^i}+S^ig_{ik}\right)
\frac \partial {\partial p_k}.
\end{equation*}
This together with the conditions (14) $b$) and (30) lead to
\begin{equation*}
g^{kj}\left( \frac{\partial \zeta _i}{\partial x^k}-\frac{\partial
\zeta _k}{
\partial x^i}+\frac{\partial S^l}{\partial y^k}g_{li}+\xi ^l\frac{\partial
g_{ik}}{\partial x^l}+S^l\frac{\partial g_{ik}}{\partial
y^l}\right) =0,
\end{equation*}
and using (31) we obtain
\begin{equation*}
\frac{\partial \zeta _i}{\partial x^k}-\frac{\partial \zeta
_k}{\partial x^i} =S(g_{ik})-2N_k^lg_{li},
\end{equation*}
and it results, considering the symmetric part (interchanging $i$
with $k$ and totalizing)
\begin{equation}
S(g_{ik})-N_k^lg_{li}-N_i^lg_{lk}=0\Leftrightarrow \nabla g\left(
\frac
\partial {\partial y^i},\frac \partial {\partial y^k}\right) =0,
\end{equation}
which means that $N$ on $TM$ is a metric nonlinear connection, and
the antisymmetric part
\begin{equation*}
\frac{\partial \zeta _i}{\partial x^k}-\frac{\partial \zeta
_k}{\partial x^i} =N_i^lg_{lk}-N_k^lg_{li},
\end{equation*}
where $g_{il}=g_{il}\circ \Psi ^{-1}=\frac{\partial \zeta
_i}{\partial y^l}$ by abuse, which is equivalent with
\begin{equation}
N_i^lg_{lk}-N_k^lg_{li}+\frac{\partial ^2L}{\partial x^k\partial
y^i}-\frac{
\partial ^2L}{\partial x^i\partial y^k}=0\Leftrightarrow \omega _L(hX,hY)=0,
\end{equation}
where $\omega _L$ is the symplectic structure induced by the
regular Lagrangian on $TM.$ From (33) and (34), using the Theorem
3.1 from \cite{Bu1} it results that $N$ is the canonical nonlinear
connection induced by the regular Lagrangian on $ TM.$ \hfill
\hbox{\rlap{$\sqcap$}$\sqcup$}

\section*{Conclusions and further developments}
The main purpose of this paper is to develop a mathematical
apparatus on the cotangent bundle as in the case of tangent bundle
and prove the dual nature of these spaces. However, on the
cotangent bundle there is a canonical symplectic form that is
independent of the hamiltonian, where as on the tangent bundle the
symplectic form is dependent to the lagrangian. Moreover, when
dealing with geodesic flows in dynamical systems, it is more
natural to think about these flows living inside the cotangent
bundle because they become automatically symplectic. \\ Thus, it
is justified to study the geometric structures on the cotangent
bundle, which are dual with the well known structures on the
tangent bundle. The existence of tensor fields $t_{ij}$ on
$T^{*}M$ (given, for example, by a pseudo-Riemannian metric or a
regular Hamiltonian) permit us to define an adapted tangent
structure $\mathcal{J}$ and a $\mathcal{J}$-regular vector field
$\rho $. A $\mathcal{J}$-regular vector field together with an
arbitrary nonlinear connection $\mathcal{N}$ define a dynamical
covariant derivative and the Jacobi endomorphism. In the case of
the horizontal $ \mathcal{J}$-regular vector field, the Jacobi
endomorphism is the contraction with $\rho $ of the curvature of
the nonlinear connection. Let us emphasize that at this point we
do not have any relation between $ (t_{ij,}\mathcal{J},\rho )$ and
the coefficients $\mathcal{N}_{ij}$ of the nonlinear connection
$\mathcal{N}.$ This will be given considering the compatibility
condition between the dynamical covariant derivative and the
adapted tangent structure, $\nabla \mathcal{J}=0$, which fix the
canonical nonlinear connection $\mathcal{N}=-\mathcal{L}_\rho
\mathcal{J}$. This canonical nonlinear connection depends only on
$t_{ij}$ and $\rho $. In this case we have the decomposition
$\nabla =\mathcal{L}_\rho +\Bbb{F}+\mathcal{J} -\Phi $ which can
be compared with the tangent case from \cite{Bu2, Ma}. Moreover,
$\nabla $ is compatible with the geometric structures induced by $
t_{ij}$ and $\rho $ (i.e. $v$, $h$, $\mathcal{J}$, $\Bbb{F}$). The
duality between the semi-Hamiltonian vector field on $T^{*}M$ and
the semispray on $ TM$ is proved using the Legendre
transformation. \\ As possible applications, we can use the
dynamical covariant derivative and Jacobi endomorphism on the
cotangent bundle in the study of symmetries of Hamitonian
dynamical systems as they are used on the tangent bundle for the
study of symmetries for systems of second order differential
equations. Also, it might be useful to find the relation between
the dynamical covariant derivative and Berwald linear connection
on the cotangent bundle.

\section*{Acknowledgments}

The author wishes to express his thanks to the referee for useful
remarks and comments concerning this paper.

Author's address:

University of Craiova,

Dept. of Statistics and Economic Informatics,

13, Al. I. Cuza, st. Craiova 200585, Romania

e-mail: liviupopescu@central.ucv.ro; liviunew@yahoo.com
\end{document}